\documentclass{article}

\usepackage{amssymb,amsmath,amsfonts}
\usepackage{url}

\newtheorem{theorem}{Theorem}[section]

\begin{document}

\title{The Tur\'an Number and Probabilistic Combinatorics}
\author{Alan J. Aw}
\date{}
\maketitle

\begin{abstract}
In this short expository article, we describe a mathematical tool called \textit{the probabilistic method}, and illustrate its elegance and beauty through proving a few well-known results. Particularly, we give an unconventional probabilistic proof of a classical theorem concerning the Tur\'{a}n number $T(n,k,l)$. 
\end{abstract}

\section{Introduction.} Probabilistic combinatorics is a field in combinatorics that was pioneered by one of the greatest mathematicians of all time, Paul Erd\H{o}s. It utilizes notions from probability theory (mainly discrete probability theory) to prove \textit{existence} results; for instance, it has provided solutions to the following questions.

\begin{itemize}
	\item For any graph $G$, can we find a bipartite subgraph which contains at least half of the edges of $G$?
	
	\item Given a family $\mathcal F$ of subsets of $[n]=\{1,2,...,n\}$, if we know its size, can we find two elements $A,B\in\mathcal F$ such that $A\subseteq B$?
\end{itemize} 

As a matter of fact, we will provide an answer to the second question above.\bigskip

Mathematicians describe the tools of probabilistic combinatorics using the generic term ``the probabilistic method.'' There is a wide range of techniques in probabilistic combinatorics, and the interested reader should consult \cite{alon}. Its diversity notwithstanding, we shall introduce the fundamental ideas of the probabilistic method; perhaps the most straightforward and unambiguous way to do so is to solve a problem. The following theorem, due to Mantel (1907) \cite{mantel}, tells us that in any graph, if it contains many edges, then it contains a triangle. Note that for illustrative purposes, we shall also present a non-probabilistic proof of the theorem.

\begin{theorem}[Mantel's Theorem]
\label{mantel}

In any graph $G(V,E)$ with $|V|=n$, if $G$ does not contain any triangle, then

\begin{equation*}
|E|\leq\left\lfloor\frac{n^2}{4}\right\rfloor.
\end{equation*}
 
\end{theorem}

\noindent\textbf{Proof (non-probabilistic).} \textit{Note: In this proof, we use $\Delta(G)$ to denote the maximum degree among the vertices of $G$.}\bigskip

Label the vertices $v_1,v_2,...,v_{n}$. Among these vertices, pick a vertex with degree $\Delta(G)$ (there may be more than one such vertex); without loss of generality, let this vertex be $v_p$ (for $1\leq p\leq n$). If $\deg(v_p)=k$, then let the set of $k$ neighbours of $v_p$ be $P$. If there exists a pair of vertices that are connected by an edge in $P$, then we are done. If not, each vertex in $P$ has degree at most $n-1-(k-1)=n-k$. Moreover, the $n-k-1$ vertices in $V\setminus P$ each have degree at most $k$. Therefore, we have

\begin{equation*}
2|E|=\sum_{i=1}^{n}{\deg(v_i)}\leq k+(n-k)k+(n-k-1)k=(2n-2k)k.
\end{equation*}

This gives us $|E|=(n-k)k\leq\left(\frac{(n-k)+k}{2}\right)^2=\frac{n^2}{4}$, where the inequality follows from the Arithmetic-Geometric Mean inequality. Since $|E|\in\mathbb N$, we have

\begin{equation*}
|E|\leq\left\lfloor\frac{n^2}{4}\right\rfloor,
\end{equation*}

\noindent as desired.\hfill $\blacksquare$\\

\noindent\textbf{Proof (probabilistic).} \textit{Note: In this proof, a clique on $k$ vertices is simply a set of $k$ vertices in which every two vertices in it are connected by an edge.}\bigskip

Consider a probability distribution $p_1,...,p_n$ on the vertex set $V$ of $G$. Let us pick two vertices $u,v$ independently and at random (allowing repetition). Then 

\begin{equation*}
\mathbb P(\{u,v\}\in E)=\sum_{i,j:\{i,j\}\in E}{p_ip_j}. 
\end{equation*}

We first assume that the probability distribution is uniform, i.e., $p_1=p_2=\cdots=p_n=\frac{1}{n}$. This gives us

\begin{equation*}
\mathbb P(\{u,v\}\in E)=\frac{2|E|}{n^2}.
\end{equation*}

Now, we modify the distribution to make $\mathbb P(\{u,v\}\in E)$ as large as possible. This happens when the probability distribution that maximizes the probability is uniform on some maximal clique. Indeed, suppose that there are two non-adjacent vertices $i,j$ such that $p_i,p_j>0$; let $s_i=\sum_{k:\{i,k\}\in E}{p_k}$ and $s_j=\sum_{k:\{j,k\}\in E}{p_k}$. If $s_i\geq s_j$ (resp. $s_i<s_j$), we set the probability of vertex $i$ to $p_i+p_j$ and the probability of vertex $j$ to zero (and conversely if $s_i<s_j$). This increases $\mathbb P(\{u,v\}\in E)$ by $s_i(p_i+p_j)-(p_is_i+p_js_j)=p_j(s_i-s_j)$ (resp. $p_i(s_j-s_i)$).\bigskip

Since the process is finite, we eventually reach a situation where there are no two non-adjacent vertices of positive probability ($p_ip_j=0$ as long as $\{i,j\}\notin E$), i.e., the probability distribution is on a clique $Q$. Then,

\begin{equation*}
\mathbb P(\{u,v\}\in E)=\mathbb P(u\neq v)=1-\sum_{i\in Q}{p_i^2}.
\end{equation*}

Assume there is no clique larger than $2$ in size (if not, then there would be a triangle violating our assumption). By Jensen's Inequality on the convex function $f(z)=z^2$, the above expression is maximized when $p_i$ is uniform on $Q$, i.e., $p_i=\frac{1}{|Q|}$, which implies

\begin{equation*}
\mathbb P(\{u,v\}\in E)\leq 1-\frac{1}{|Q|}\leq 1-\frac{1}{2},
\end{equation*}

assuming that there is no clique larger than $2$ in size (if not there would be a triangle, violating our assumption). Since we started with $\mathbb P(\{u,v\}\in E)=\frac{2|E|}{n^2}$ and never decreased it in the process, we yield

\begin{equation*}
|E|\leq\frac{n^2}{4}\Longrightarrow |E|\leq\left\lfloor\frac{n^2}{4}\right\rfloor,
\end{equation*}

\noindent as desired.\hfill $\blacksquare$\\ 

In fact, Theorem \ref{mantel} can be generalized to the same problem for cliques of size $k$, where $k\geq3$ (for that matter, Mantel established the result for $k=2$). The general problem was solved by the Hungarian mathematician Paul Tur\'{a}n, and a wonderful probabilistic proof of it can be found in \cite{aigner}.

\section{Further Examples.} Recall the second question in the Introduction regarding the size of a family of sets of $[n]$. Intuitively, if the size of the family, $|\mathcal F|$, is large enough, then we should be able to find some $A,B\in\mathcal F$ such that $A\subset B$. Indeed, the following theorem of Sperner gives us a sufficient condition on $\mathcal F$ such that we can find two sets such that one contains the other.

\begin{theorem}[Sperner's Theorem]
\label{sperner}

If $\mathcal F$ is a Sperner family of subsets (or an antichain) of $[n]$ (i.e., for all $A,B\in\mathcal F$ we have $A\not\subset B$ and $B\not\subset A$), then

\begin{equation*}
|\mathcal F|\leq{n\choose\lfloor\frac{n}{2}\rfloor}.
\end{equation*}

\end{theorem} 

\noindent\textbf{Proof.} Without loss of generality, let $\mathcal F=\{A_1,A_2,...,A_{|\mathcal F|}\}$. Consider any $A_i\in\mathcal F$ with $|A_i|=a_i$. Randomly and uniformly pick any permutation of $[n]$. Denoting by $E_i$ the event that the elements in $A_i$ appear in the first $a_i$ numbers in the permutation, we get

\begin{equation*}
\mathbb P(E_i)=\frac{(a_i)!(n-a_i)!}{n!}=\frac{1}{{n\choose a_i}}\geq\frac{1}{{n\choose{\lfloor\frac{n}{2}\rfloor}}}.
\end{equation*}

Moreover, the $E_i$ are mutually exclusive because $E_i\cap E_j$ is the event that both $A_i$ and $A_j$ appear at first in the permutation, implying that one of $A_i\subset A_j$ or $A_j\subset A_i$ holds. This contradicts our hypothesis. Therefore, summing up our probabilities and noting that probabilities are always less than 1, we obtain

\begin{equation*}
1\geq\mathbb P\left(\bigcup_{i=1}^{|\mathcal F|}{E_i}\right)=\sum_{i=1}^{|\mathcal F|}{\frac{1}{{n\choose a_i}}}\geq \frac{|\mathcal F|}{{n\choose {\lfloor\frac{n}{2}\rfloor}}},
\end{equation*}

\noindent where the last inequality follows from the fact that ${n\choose p}$ is maximized when $p=\left\lfloor\frac{n}{2}\right\rfloor$. The result follows.\hfill $\blacksquare$\\

Next, we have a classical result that is related to the Tur\'{a}n number $T(n,k,l)$. We shall present a probabilistic proof which differs from the conventional non-probabilistic proof (see for instance \cite[chapter 1]{jukna}) of the theorem.

\begin{theorem}
\label{turan}

The Tur\'{a}n number $T(n,k,l)$, where $n\geq k\geq l$, is the smallest number of $l$-element subsets of an $n$-element set $X$ such that every $k$-element subset of $X$ contains at least one of these sets. If $T(n,k,l)=|\mathcal F|$, where $\mathcal F$ is the smallest family satisfying the above conditions, then

\begin{equation*}
|\mathcal F|\geq\frac{{n\choose l}}{{k\choose l}}.
\end{equation*}

\end{theorem}

\noindent\textbf{Proof.} Let us randomly select without replacement a $k$-element subset $K$ from $X$. Let $A$ be the event that the first $l$ elements picked form a set $L$ which belongs to $\mathcal F$. Then

\begin{equation*}
\mathbb P(A)=|\mathcal F|\cdot\left(\frac{l}{n}\cdot\frac{l-1}{n-1}\cdots\frac{1}{n-l+1}\cdot\frac{n-l}{n-l}\cdots\frac{n-k+1}{n-k+1}\right)=\frac{|\mathcal F|}{{n\choose l}}.
\end{equation*}

On the other hand, let us consider the events $E_i$ where $E_i$ is the event that the $i$th $k$-element subset was selected. Then we have

\begin{equation*}
\mathbb P(A)=\sum_{i=1}^{{n\choose k}}{\mathbb P(A\cap E_i)}\geq\sum_{i=1}^{{n\choose k}}{\left[\frac{1}{{n\choose k}}\cdot\frac{1}{{k\choose l}}\right]}=\frac{1}{{k\choose l}},
\end{equation*}

\noindent since each $k$-element subset of $X$ contains at least one element in $\mathcal F$. Combining the two equations yields the desired inequality.\hfill $\blacksquare$

\section{Conclusion} Truly, the probabilistic method presents a pristine and elegant approach to problem solving. From its very first applications in extremal graph theory by Paul Erd\H{o}s, to its plethora of applications (see \cite{alon,jukna}) in coding theory, number theory and geometry, probabilistic combinatorics today has certainly emerged as one of the richest subdisciplines of mathematics where beauty, creativity and rigour converge. In Erd\H{o}s' own words, ``This one's from The Book!''

\paragraph{Acknowledgments.} The author would like to thank Professor Sergei Tabachnikov for his comments that helped make several improvements to this article.

\bigskip

\noindent\textit{Raffles Institution,
One Raffles Institution Lane, S575954,\\
Singapore\\
nalawanij@gmail.com}

\end{document}